\documentclass{amsart}
\usepackage[all]{xy}
\usepackage{comment}

\newtheorem{theorem}{Theorem}[section]
\newtheorem{lemma}{Lemma}[section]
\newtheorem{proposition}{Proposition}[section]

\newtheorem{definition}{Definition}

\theoremstyle{remark}
\newtheorem{remark}{Remark}[section]
\numberwithin{equation}{section}
\numberwithin{lemma}{section}
\numberwithin{proposition}{section}
\numberwithin{definition}{section}

\newcommand{\internalcomment}[1]{}

\begin{document}

\title[Homotopy theory of triangulated categories]{Homotopy theory of
  well-generated algebraic triangulated categories}
\author{Gon{\c c}alo Tabuada}
\address{Universit{\'e} Paris 7 - Denis Diderot, UMR 7586
  du CNRS, case 7012, 2 Place Jussieu, 75251 Paris cedex 05, France.}

\thanks{Supported by FCT-Portugal, scholarship {\tt SFRH/BD/14035/2003}.}
\keywords{Monads and algebras, Quillen model structure, well-generated algebraic
  triangulated category, Verdier's quotient, exact DG category, $\alpha$-cocompletness}

\email{
\begin{minipage}[t]{5cm}
tabuada@math.jussieu.fr
\end{minipage}
}

\begin{abstract}
For every regular cardinal $\alpha$, we construct a cofibrantly
generated Quillen model structure on a category whose objects are essentially DG
categories which are stable under suspensions, cosuspensions, cones
and $\alpha$-small sums.

Using results of Porta, we show that the category
of well-generated (algebraic) triangulated categories in the sense of
Neeman is naturally enhanced by our Quillen model category.
\end{abstract}

\maketitle

\tableofcontents

\section{Introduction}
Triangulated categories appear naturally in several branches
of mathematics such as algebraic geometry, algebraic analysis,
$K$-theory, representation theory and even in mathematical physics,
see for instance \cite{Kontsevich}.

In his book~\cite{Neeman}, Neeman introduces an important class of
triangulated categories known as {\it well-generated}. Recall that a
triangulated category $\mathcal{T}$ is well-generated if it is
$\alpha$-compactly generated for a regular cardinal $\alpha$,
see~\cite{Krause} \cite{Neeman}. Neeman proves that the Brown representability
theorem holds for well-generated triangulated categories and that this
class of triangulated categories is stable under localizations and
passage to localizing subcategories generated by a set of objects.

In this paper we study the category of well-generated algebraic triangulated categories, see \cite{dg-cat-survey}, using the tools of Quillen's homotopical
algebra, \cite{Quillen}, and the formalism of DG categories,
\cite{Drinfeld}~\cite{dg-cat-survey}~\cite{ICM}~\cite{Toen}.
More precisely, for a fixed regular cardinal $\alpha$ we construct a
category $\mathsf{dgcat}_{ex,\alpha}$, whose objects are essentially
the dg categories which are stable under suspensions, cosuspensions,
cones and $\alpha$-small sums.

The construction of $\mathsf{dgcat}_{ex,\alpha}$ is done in two steps.

First we construct a monad $\mathsf{T}_{\alpha}$ on the
  category of small differential graded categories $\mathsf{dgcat}$
  and we consider the associated category $\mathsf{T}_{\alpha}$-$\mathsf{alg}$ of
  $\mathsf{T}_{\alpha}$-algebras.

Then the category $\mathsf{dgcat}_{ex,\alpha}$, is obtained by
  considering specific diagrams in $\mathsf{T}_{\alpha}$-$\mathsf{alg}$. 

In each of these two steps we dispose of an adjunction and by applying
an argument due to Quillen to our situation, we are able to lift the Quillen
model structure on $\mathsf{dgcat}$, see~\cite{cras}, along these
adjunctions to $\mathsf{dgcat}_{ex,\alpha}$.

Finally, we define a functor $\mathcal{D}_{\alpha}(-)$ from
$\mathsf{dgcat}_{ex, \alpha}$ to the category $\mathsf{Tri}_{\alpha}$
of $\alpha$-compactly generated triangulated categories,
see~\cite{Krause} \cite{Neeman}, which by \cite{Ober}~\cite{Porta}
is known to verify the following conditions:
\begin{enumerate}
\item[-] every category $\mathcal{T}$ in $\mathsf{Tri}_{\alpha}$ is equivalent to
  $\mathcal{D}_{\alpha}(\underline{A})$ for some $\underline{A}$ in $\mathsf{dgcat}_{ex,\alpha}$ and
\item[-] a morphism $F$ in $\mathsf{dgcat}_{ex, \alpha}$ is a
  weak equivalence if and only if $\mathcal{D}_{\alpha}(F)$ is an
  equivalence of triangulated categories.
\end{enumerate}

This shows that well-generated algebraic triangulated categories up to
equivalence admit a natural Quillen enhancement given by our model
category $\mathsf{dgcat}_{ex, \alpha}$.

\section{Acknowledgments}
This article is part of my Ph.D. thesis under the supervision of
Prof. B.~Keller. I deeply thank him for suggesting this problem to me and for
several useful discussions. I am also very grateful to M.~Batanin and
M.~Weber for e-mail exchange concerning monadic structures and to
C.~Berger and B.~To{\"e}n for useful discussions concerning homotopical algebra.

\section{Preliminaries}
In what follows $\alpha$ will denote a regular cardinal and $k$ a
commutative ring with unit. The tensor product $\otimes$ will denote the tensor product over $k$. Let
$\mathsf{Ch}(k)$ denote the category of complexes over $k$. By a
{\it dg category}, we mean a differential graded $k$ category, see
\cite{Drinfeld} \cite{ICM} \cite{cras}. For a dg category
$\mathcal{A}$, we denote by $\mathcal{C}_{dg}(\mathcal{A})$ the dg
category of right $\mathcal{A}$ dg modules and by $\,\widehat{ } : \mathcal{A} \rightarrow
\mathcal{C}_{dg}(\mathcal{A})$ the Yoneda dg functor.
We write $\mathsf{dgcat}$ for the category of small dg categories.
It is proven in \cite{cras}, that the
category $\mathsf{dgcat}$ admits a structure of cofibrantly generated
model category whose weak equivalences are the quasi-equivalences.

\section{Monadic structure $\mathsf{T}$ on $\mathsf{dgcat}$}\label{monad}

In this section, we will construct a monad $\mathsf{T}_{\alpha}$ on
$\mathsf{dgcat}$. The $\mathsf{T}_{\alpha}$-algebras will be
essentially the dg categories which admit $\alpha$-small
sums, see proposition~\ref{sums}. This monad $\mathsf{T}_{\alpha}$ is
a variant of the well-known coproduct complection monad $\mathsf{Fam}$
of `families', \emph{cf.}~\cite{Fam}, see also~\cite{Kelly}.

For each ordinal $\beta$ strictly smaller than $\alpha$, we denote by
$I_{\beta}$ the underlying set of $\beta$, see
\cite{Set}~\cite{Hirschhorn}.
Let $\mathcal{A}$ be a small dg category.

\begin{definition}\label{monadf}
Let $\mathsf{T}_{\alpha}(\mathcal{A})$ be the small dg category whose
objects are the maps
$$ I_{\beta} \stackrel{F_{\beta}}{\longrightarrow}
\mbox{obj}\,(\mathcal{A})\,,$$
where $\beta < \alpha$, and whose complex of morphisms between $F_{\beta}$ and $F_{\beta '}$
is 
$$ \mathsf{Hom}_{\mathsf{T}_{\alpha}(\mathcal{A})}(F_{\beta}, F_{\beta
  '}):= \prod_{i \in I_{\beta}}  \bigoplus_{j \in I_{\beta'}} \mathsf{Hom}_{\mathcal{A}}(F_{\beta}(i), F_{\beta
'}(j))\,.$$
The composition and the identities are induced by those of $\mathcal{A}$.
\end{definition}

\begin{remark}
Clearly $\mathsf{T}_{\alpha}(\mathcal{A})$ is a small dg category and
the above construction is functorial in $\mathcal{A}$. We dispose of a
functor
$$ \mathsf{T}_{\alpha}(-): \mathsf{dgcat} \longrightarrow \mathsf{dgcat}\,.$$
\end{remark}

\begin{definition}
Let $\eta$ be the natural transformation
$$ \eta: Id_{\mathsf{dgcat}} \longrightarrow
\mathsf{T}_{\alpha}(-)\,,$$
whose evaluation at $\mathcal{A}$ is the dg functor
$$ \eta_{\mathcal{A}}:\mathcal{A} \longrightarrow
\mathsf{T}_{\alpha}(\mathcal{A})$$
that sends an object $X$ of $\mathcal{A}$ to the map
$$ I_{\mathbf{1}} \stackrel{X}{\longrightarrow}{
  \mbox{obj}\,(\mathcal{A})}\,,$$
where $\mathbf{1}$ denotes the first sucessor ordinal, see \cite{Set}.
\end{definition}
Remark that $\eta_{\mathcal{A}}$ is a fully faithful dg functor.

\begin{definition}\label{assoc}
Let $\mu_{\mathcal{A}}$ be the dg functor
$$ \mu_{\alpha}:(\mathsf{T}_{\alpha}\circ
\mathsf{T}_{\alpha})(\mathcal{A}) \longrightarrow
\mathsf{T}_{\alpha}(\mathcal{A})\,,$$
which sends an object
$$
\begin{array}{ccc}
I_{\beta} & \stackrel{F_{\beta}}{\longrightarrow} &
\mbox{obj}\,(\mathsf{T}_{\alpha}(\mathcal{A})) \\
x & \longmapsto & (I_{\gamma_x}
\stackrel{F_{\gamma_x}}{\longrightarrow} \mbox{obj}\,(\mathcal{A}))\,,
\end{array}
$$
of $(\mathsf{T}_{\alpha} \circ \mathsf{T}_{\alpha})(\mathcal{A})$ to
the map
$$ \coprod_{x \in \beta} F_{\gamma_x} := I_{\underset{x \in \beta}{\sum}
  \gamma_x} \longrightarrow \mbox{obj} \,(\mathcal{A})\,,$$
where $\underset{x \in \beta}{\sum} \gamma_x$ denotes the increasing sum of the
ordinals $\gamma_x$, $x \in \beta$, along the ordinal $\beta$, see
\cite{Set}.
\end{definition}

\begin{remark}
Remark that the ordinal $\underset{x \in \beta}{\sum}\gamma_x$ is strictly
smaller than $\alpha$ since $\beta$ and each one of the $\gamma_x$'s
are strictly smaller than $\alpha$, which is by hypothesis a regular
cardinal. The above construction is functorial in $\mathcal{A}$ and so we dispose of
a natural transformation
$$\mu:(\mathsf{T}_{\alpha} \circ \mathsf{T}_{\alpha})(-) \longrightarrow \mathsf{T}_{\alpha}(-)\,.$$
\end{remark}

\begin{proposition}
We dispose of a monad, see~\cite{Macl},
$\mathsf{T}_{\alpha}=(\mathsf{T}_{\alpha}(-), \eta, \mu)$ on the
category $\mathsf{dgcat}$.
\end{proposition}

\begin{proof}
We need to prove that the diagrams of functors
$$
\xymatrix{
\mathsf{T}_{\alpha}\circ \mathsf{T}_{\alpha}\circ \mathsf{T}_{\alpha}
\ar[d]_{\mu \mathsf{T}_{\alpha}} \ar[r]^-{\mathsf{T}_{\alpha}\mu} &
\mathsf{T}_{\alpha}\circ \mathsf{T}_{\alpha} \ar[d]^{\mu} & &
\mathsf{T}_{\alpha} \ar[r]^-{\eta \mathsf{T}_{\alpha}} \ar@{=}[dr] & \mathsf{T}_{\alpha} \circ
\mathsf{T}_{\alpha} \ar[d]^{\mu} & \mathsf{T}_{\alpha} \ar[l]_-{\mathsf{T}_{\alpha} \eta}
\ar@{=}[dl] \\
\mathsf{T}_{\alpha}\circ \mathsf{T}_{\alpha} \ar[r]_{\mu} & \mathsf{T}_{\alpha} & & &
\mathsf{T}_{\alpha} & 
}
$$
are commutative. The left diagram is commutative because the ordinal
sum operation, see~\cite{Set}, on the set of ordinals strictly smaller
than $\alpha$, is associative. The left triangle of the right diagram
is also commutative because the ordinal sum operation admits a
unity which is given by the first ordinal. The right triangle of the right
diagram is also commutative by definition of $\mu$.

This proves the proposition.
\end{proof}

Let $\mathcal{A}$ be a small dg category.

\begin{lemma}\label{sums1}
The dg category $\mathsf{T}_{\alpha}(\mathcal{A})$ admits
$\alpha$-small sums.
\end{lemma}

\begin{proof}
Let $\beta$ be a cardinal strictly smaller than $\alpha$. Let
$J$ be a set of cardinality $\beta$ and $G$ a morphism
$$ G: J \longrightarrow
\mbox{obj}\,(\mathsf{T}_{\alpha}(\mathcal{A}))\,.$$
Choose a bijection between $J$ and $I_{\beta}$ and consider the
object 
$$ F_{\beta}:I_{\beta} \longrightarrow
\mbox{obj}\,(\mathsf{T}_{\alpha}(\mathcal{A}))$$
of $(\mathsf{T}_{\alpha} \circ \mathsf{T}_{\alpha})$ associated to $G$. Now remark that
by definition of $\mathsf{T}_{\alpha}(\mathcal{A})$ and since 
$$I_{\underset{x \in \beta}{\sum}\gamma_x}= \underset{x \in
  I_{\beta}}{\cup} I_{\gamma_x} $$
the object $\mu_{\mathcal{A}}(F_{\beta})$ of
$\mathsf{T}_{\alpha}(\mathcal{A})$ is in fact the $\beta$-small sum of $G$.
\end{proof}

Recall from \cite{Macl} that by definition an algebra over this monad (=$\mathsf{T}_{\alpha}$-algebra) consists of a pair $A=(\mathcal{A},R)$, where
$\mathcal{A}$ is a small dg category and $R$ is a dg functor
$R:\mathsf{T}_{\alpha}(\mathcal{A}) \rightarrow \mathcal{A}$ which makes both diagrams
$$
\xymatrix{
(\mathsf{T}_{\alpha}\circ \mathsf{T}_{\alpha})(\mathcal{A})
\ar[d]_{\mu_{\mathcal{A}}} \ar[r]^-{\mathsf{T}_{\alpha}(R)} &
\mathsf{T}_{\alpha}(\mathcal{A}) \ar[d]^{R} & &
\mathcal{A} \ar@{=}[dr] \ar[r]^-{\eta_{\mathcal{A}}} & \mathsf{T}_{\alpha}(\mathcal{A}) \ar[d]^{R} \\
\mathsf{T}_{\alpha}(\mathcal{A}) \ar[r]_{R} & \mathcal{A} & & &
\mathcal{A} 
}
$$
commute. A morphism $F:(\mathcal{A},R) \rightarrow (\mathcal{B},G)$ of
$\mathsf{T}_{\alpha}$-algebras is a dg functor $F:\mathcal{A}
\rightarrow \mathcal{B}$ which renders commutative the diagram
$$
\xymatrix{
\mathsf{T}_{\alpha}(\mathcal{A}) \ar[d]_R
\ar[r]^{\mathsf{T}_{\alpha}(F)} & \mathsf{T}_{\alpha}(\mathcal{B})
\ar[d]^G \\
\mathcal{A} \ar[r]_F & \mathcal{B}\,.
}
$$
We denote by $\mathsf{T}_{\alpha}$-$\mathsf{alg}$ the category of
$\mathsf{T}_{\alpha}$-algebras. By theorem $1$ of chapter VI from \cite{Macl},
we dispose of an adjunction
$$
\xymatrix{
\mathsf{T}_{\alpha}\mbox{-}\mathsf{alg} \ar@<1ex>[d]^U \\
\mathsf{dgcat} \ar@<1ex>[u]^F\,,
}
$$
where the functor $U$ associates to a $\mathsf{T}_{\alpha}$-algebra $A$ the dg category $\mathcal{A}$ and $F$ associates to a dg category $\mathcal{B}$ the
$\mathsf{T}_{\alpha}$-algebra $(\mathsf{T}_{\alpha}(\mathcal{B}),\mu_{\mathcal{B}})$. 

\begin{proposition}\label{sums}
Let $A=(\mathcal{A},R)$ be a $\mathsf{T}_{\alpha}$-algebra. Then the dg
category $\mathcal{A}$ admits $\alpha$-small sums.
\end{proposition}

\begin{proof}
Let $\beta$ be a cardinal strictly smaller than $\alpha$ and $J$
a set of cardinality $\beta$. Let $G$ be a morphism
$$ G:J \longrightarrow \mbox{obj}\,(\mathcal{A})\,.$$
Choose a bijection between $J$ and $I_{\beta}$ and consider the
object of $\mathsf{T}_{\alpha}(\mathcal{A})$
$$ F_{\beta}:I_{\beta} \longrightarrow \mbox{obj}\,(\mathcal{A})\,,$$ associated to $G$.
We will prove that the object $R(F_{\beta})$ of $\mathcal{A}$ is the
$\beta$-small sum of $G$. Consider the following object 
$$ \overline{F_{\beta}}:=
\mathsf{T}_{\alpha}(\eta_{\mathcal{A}})(F_{\beta})\,,$$
in $(\mathsf{T}_{\alpha}\circ \mathsf{T}_{\alpha})(\mathcal{A})$.
Remark that $F_{\beta} = \mu_{\mathcal{A}}(\overline{F_{\beta}})$ and so
  by lemma~\ref{sums1}, $F_{\beta}$ is the $\beta$-small sum of
  $\overline{F_{\beta}}$. Now, since $A$ is a $\mathsf{T}_{\alpha}$-algebra the
  equality $R \circ \eta_{\mathcal{A}} = Id$ implies that
  $R(F_{\beta})$ is a weak $\beta$-small sum of $G$ in
  $\mathcal{A}$. We now prove that this weak $\beta$-small sum is in
  fact a true one. Remark that for each $x \in I_{\beta}$, we dispose of a
  closed morphism of degree zero
$$ i_x: F_{\beta}(x) \longrightarrow R(F_{\beta})\,,$$
in $\mathcal{A}$.
Now let $Z$ be an object of $\mathcal{A}$ and $H_1$ and $H_2$ two
elements of $\mathsf{Hom}_{\mathcal{A}}(R(F_{\beta}),Z)$ such that
$$H_1 \circ i_x = H_2 \circ i_x, \, \forall x \in I_{\beta}\,.$$
We will now prove that $H_1=H_2$.

Consider the following commutative diagram in $\mathcal{A}$
$$
\xymatrix{
F_{\beta}(x) \ar[r]^-{R(i_x)} \ar[dr]_{G_x} & R(F_{\beta}) \ar[d]_{H_1} \ar@<1ex>[d]^{H_2}\\
 & Z
}
$$
Apply the dg functor $\eta_{\mathcal{A}}$ to the previous diagram and
consider the following one
$$
\xymatrix{
 &  & F_{\beta} \ar@{.>}[d]^-{\theta} \ar@/^5pc/@{.>}[dd]^-{\Phi}\\
\eta_{\mathcal{A}}(F_{\beta}(x)) \ar[urr]^-{i_x}
\ar[rr]^-{\eta_{\mathcal{A}}(R(i_x))}
\ar[drr]_-{\eta_{\mathcal{A}}(G_x)} & &
\eta_{\mathcal{A}}(R(F_{\beta})) \ar[d]_-{\eta_{\mathcal{A}}(H_1)}
\ar@<1ex>[d]^-{\eta_{\mathcal{A}}(H_2)} \\
 & & \eta_{\mathcal{A}}(Z) \,,
}
$$
where $\theta$ is the unique morphism in
$\mathsf{T}_{\alpha}(\mathcal{A})$ such that
$$ i_x \circ \theta = \eta_{\mathcal{A}}(R(i_x)), \, \forall x \in
I_{\beta}$$
and $\Phi$ is the unique morphism in
$\mathsf{T}_{\alpha}(\mathcal{A})$ such that 
$$ i_x \circ \Phi = \eta_{\mathcal{A}}(G_x), \, \forall x \in
I_{\beta}\,.$$
This implies that
$$ \eta_{\mathcal{A}}(H_1)\circ \theta = \eta_{\mathcal{A}}(H_2)\circ
\theta \,.$$
We will show that $R(\theta)=Id$, which immediately implies the
proposition. Recall that since $A$ is a $\mathsf{T}_{\alpha}$-algebra, we
dispose of the following commutative diagram
$$
\xymatrix{
(\mathsf{T}_{\alpha}\circ \mathsf{T}_{\alpha})(\mathcal{A})
\ar[d]_{\mu_{\mathcal{A}}} \ar[r]^-{\mathsf{T}_{\alpha}(R)} &
\mathsf{T}_{\alpha}(\mathcal{A}) \ar[d]^{R} \\
\mathsf{T}_{\alpha}(\mathcal{A}) \ar[r]_{R} & \mathcal{A}\,.
}
$$
Remark that since  $\overline{F_{\beta}}$ is the $\beta$-small sum of the objects
$\eta_{\mathsf{T}_{\alpha}(\mathcal{A})}(\eta_{\mathcal{A}}(F_{\beta}(x)))$,
the morphisms $i_x$ induce a morphism
$$ \Psi: \overline{F_{\beta}} \longrightarrow
\eta_{\mathsf{T}_{\alpha}(\mathcal{A})}(F_{\beta})\,,$$
in $(\mathsf{T}_{\alpha} \circ \mathsf{T}_{\alpha})(\mathcal{A})$.
Finally remark that $\mathsf{T}_{\alpha}(R)(\Psi)=\theta$ and $\mu_{\mathcal{A}}(\Psi)=Id$. This implies that $R(\theta)=Id$ and so
the proposition is proven.
\end{proof}

\begin{remark}\label{sums2}
Remark that the proof of the proposition~\ref{sums} shows us that if
$F: A \rightarrow B$ is a morphism of $\mathsf{T}_{\alpha}$-algebras, then the dg functor $F:\mathcal{A} \rightarrow \mathcal{B}$ preserves
$\alpha$-small sums.
\end{remark}

\section{Quillen's lifting argument}

In this section, we consider an argument due to Quillen,
see~\cite{Quillen}, that allows us to lift a Quillen model structure
along an adjunction.

Let $\mathcal{N}$ be a complete and cocomplete category. Consider a functor
$U:\mathcal{N} \rightarrow \mathcal{M}$, with $\mathcal{M}$ a Quillen
model category. Assume that $U$ admits a left adjoint 
$$F:\mathcal{M}\rightarrow \mathcal{N}\,.$$

\begin{definition}\label{lift}
A morphism $f:A \longrightarrow B$ in $\mathcal{N}$ is:
\begin{enumerate}
\item[-] a weak equivalence if $U(f)$ is a quasi-equivalence in
  $\mathcal{M}$.
\item[-] a fibration if $U(f)$ is a fibration in $\mathcal{M}$.
\item[-] a cofibration if it has the left lifting property with respect
  to all trivial fibrations in $\mathcal{N}$.
\end{enumerate}
\end{definition}

\begin{theorem}\label{liftarg}
Suppose that $\mathcal{M}$ is a cofibrantly generated Quillen model
category and that $U$ commutes with $\alpha$-filtered colimits for a
regular cardinal $\alpha$, see~\cite{Hirschhorn}.

Then with the notions of weak equivalence, fibration and cofibration
defined above, $\mathcal{N}$ is a Quillen model category provided the
following assumption on cofibrations holds: every cofibration with the
left lifting property with respect to fibrations is a weak equivalence.
\end{theorem}

\begin{proof}
We denote by $I$ the set of generating cofibrations of $\mathcal{M}$
and by $J$ the set of generating trivial cofibrations of
$\mathcal{M}$. Remark that since $\mathcal{M}$ is cofibrantly
generated, the domains of the elements of the sets $I$ and $J$ are
$\beta$-small for a cardinal $\beta$, see \cite{Hirschhorn}. 
Since $U$ commutes with $\alpha$-small filtered colimits, the images of
these domains under the functor $F$ will be $\gamma$-small, where
$\gamma$ is the maximum of $\alpha$ and $\beta$.
This shows that the sets $F(I)$ and $F(J)$ of $\mathcal{N}$ allow the
small object argument.

Now the proof follows the lines of the one of
theorem $4.1$ from \cite{Jardine}: We simply use the set $I$, respectively $J$, instead
of the generating cofibrations of simplicial sets, respectively generating acyclic
cofibrations of simplicial sets and consider $\gamma$-transfinite
compositions, see~\cite{Hirschhorn}, for the construction of the
factorizations.

Remark also that the class of cofibrations that have the left lifting
property with respect to fibrations is stable under $\gamma$-transfinite
compositions.
\end{proof}

Now let $\mathcal{N}$ and $\mathcal{M}$ be as at the beginning of this
section and consider definition~\ref{lift}.
 
\begin{proposition}\label{Quillen}
Suppose that
\begin{itemize}
\item[-] for every object $A$ in $\mathcal{N}$, the unique morphism $A
  \rightarrow *$, where $*$ denotes the terminal object in
  $\mathcal{N}$, is a fibration.
\item[-] for every object $A$ in $\mathcal{N}$, we dispose of a
  factorization
$$
\xymatrix{
A \ar[rr]^{\Delta} \ar[dr]^{\sim}_{i_A} & & A \times A \\
 & P(A) \ar@{->>}[ur]_{q_A} & \,,
}
$$
where $i_A$ is a weak equivalence and $q_A$ is a fibration.
\end{itemize}
Then every morphism in $\mathcal{N}$ that has the left lifting property with
respect to all fibrations is a weak equivalence.
\end{proposition}

\begin{proof}
Let $i:A \rightarrow B$ be a morphism in $\mathcal{N}$ that
has the left lifting property with respect to all fibrations.
Consider the following diagram
$$
\xymatrix{
A \ar@{=}[r] \ar[d]_i &  A \ar@{->>}[d]\\
B \ar@{->>}[r] \ar@{.>}[ur]^u & \ast \,.
}
$$
By the hypothesis on $i$ we dispose of a morphism $u$ such that $u\circ
i=\mbox{Id}$. We now show that the morphism $i \circ u$ is right
homotopic to the
identity of $B$.
By hypothesis, we dispose of a factorization
$$
\xymatrix{
B \ar[rr]^{\Delta} \ar[dr]^{\sim}_{i_B} & & B\times B \\
 & P(\mathcal{B}) \ar@{->>}[ur]_{q_B} & 
}
$$
which allows us to construct the diagram
$$
\xymatrix{
A \ar[rr]^{i_B \circ i} \ar[d]_i && P(\mathcal{B}) \ar@{->>}[d]^{q_B} \\
B \ar[rr]_{[Id\, , \, i\circ u]} \ar@{.>}[urr]^H & & B\times B\,.
}
$$
By the hypothesis on $i$, we dispose of a morphism $H$, which by
definition is a right homotopy between the identity of $B$ and
$i\circ u$. Since the functor $U: \mathcal{N} \rightarrow \mathcal{M}$ preserves products, fibrations and weak
equivalences the identity on  $U(B)$ and $U(i)\circ U(u)$
are right homotopic in $\mathcal{M}$ and so they become equal in the
homotopy category $\mathsf{Ho}(\mathcal{M})$. Since we already now that
$U(u)\circ U(i)$ is the identity on $U(B)$, we conclude that the morphism $U(i)$ is an
isomorphism in $\mathsf{Ho}(\mathcal{M})$. By proposition $1.14$ from
\cite{Jardine}, $U(i)$ is in fact a weak equivalence in
$\mathcal{M}$ which implies by definition that $i$ is a weak
equivalence in $\mathcal{N}$. This proves the lemma.
\end{proof}

\section{Homotopy theory of $\mathsf{T}$-algebras}\label{homotopy}

Recall from section~\ref{monad} that we dispose of an adjunction
$$
\xymatrix{
\mathsf{T}_{\alpha}\mbox{-}\mathsf{alg} \ar@<1ex>[d]^U \\
\mathsf{dgcat} \ar@<1ex>[u]^F\,.
}
$$
Since the category $\mathsf{dgcat}$ is complete, proposition $4.3.1$
from \cite{Bor} implies that $\mathsf{T}_{\alpha}$-$\mathsf{alg}$ is
also complete. Now remark that the functor
$$ \mathsf{T}_{\alpha}(-): \mathsf{dgcat} \longrightarrow
\mathsf{dgcat}\,,$$
see definition~\ref{monadf}, commutes with $\alpha$-filtered
colimits, see \cite{Bor}. This implies by proposition $4.3.2$ and
$4.3.6$ from \cite{Bor} that the category $\mathsf{T}_{\alpha}$-$\mathsf{alg}$ is
cocomplete and that the functor $U$ commutes with $\alpha$-filtered
colimits.

From now on and until the end of this section we consider the
definition~\ref{lift} applied to our particular adjunction $(F,U)$.

Let $\mathcal{B}$ be a small dg category.

\begin{definition}
Let $P(\mathcal{B})$ be the dg category, see~\cite{Drinfeld}, whose
objects are the closed morphismes of degree zero in $\mathcal{B}$
$$ X \stackrel{f}{\longrightarrow} Y\,,$$
that become invertible in $\mathsf{H}^0(\mathcal{B})$.
We define the complex of morphismes
$$ \mathsf{Hom}_{P(\mathcal{B})}(X\stackrel{f}{\rightarrow}Y,
X'\stackrel{f'}{\rightarrow}Y')$$
as the homotopy pull-back in $\mathsf{Ch}(k)$ of the diagram
$$
\xymatrix{
& \mathsf{Hom}_{\mathcal{B}}(Y,Y') \ar[d]^{f^*} \\
\mathsf{Hom}_{\mathcal{B}}(X,X') \ar[r]^{f'_*} & \mathsf{Hom}_{\mathcal{B}}(X,Y')\,.
}
$$
\end{definition}

It is proven in \cite{dgquot} that $P(\mathcal{B})$ is a path object
for $\mathcal{B}$ in the Quillen model structure on $\mathsf{dgcat}$,
see~\cite{cras}. Remark that the above construction is functorial in
$\mathcal{B}$ and so we dispose of a functor
$$P(-):\mathsf{dgcat} \rightarrow \mathsf{dgcat}\,.$$

Let $B=(\mathcal{B},S)$ be a $\mathsf{T}_{\alpha}$-algebra.

\begin{proposition}\label{path}
The category $P(\mathcal{B})$ carries a natural
$\mathsf{T}_{\alpha}$-algebra structure and so $B$ admits a
path-object $P(B)$.
\end{proposition}

\begin{proof}
We dispose of a dg functor
$$ S: \mathsf{T}_{\alpha}(\mathcal{B}) \longrightarrow \mathcal{B}$$
and we will construct a dg functor $\overline{S}$ from
$\mathsf{T}_{\alpha}(P(\mathcal{B})) \longrightarrow
P(\mathcal{B})\,.$
Remark that we dispose of a faithful dg functor
$$ \mathsf{T}_{\alpha}(P(\mathcal{B})) \longrightarrow P(\mathsf{T}_{\alpha}(\mathcal{B}))\,.$$
Define $\overline{S}$ as the composition of this dg functor with 
$$ P(\mathsf{T}_{\alpha}(\mathcal{B})) \longrightarrow
P(\mathcal{B})\,.$$
Since $B$ is a $\mathsf{T}_{\alpha}$-algebra, this construction shows us that
$(P(\mathcal{B}), \overline{S})$ is also a $\mathsf{T}_{\alpha}$-algebra. It is
also clear by construction that the dg functors $\mathcal{A}
\stackrel{i_{\mathcal{A}}}{\longrightarrow} P(\mathcal{A})$ and 
$$
\xymatrix{
 P(\mathcal{A}) \ar@{->>}[r]^-{q_{\mathcal{A}}} &
 \mathcal{A}\times\mathcal{A}
}
$$ are in fact morphismes of
$\mathsf{T}_{\alpha}$-algebras, where $\mathcal{A}\times\mathcal{A}$
carries the diagonal $\mathsf{T}_{\alpha}$-action. This proves the proposition. 
\end{proof}

\begin{theorem}\label{main}
The category $\mathsf{T}_{\alpha}$-$\mathsf{alg}$ when endowed with
the notions of weak equivalence, fibration and cofibration as in
definition~\ref{lift}, becomes a cofibrantly generated Quillen model
category and the adjunction $(F,U)$ becomes a Quillen adjunction.
\end{theorem}

\begin{proof}
Recall from \cite{cras} that we dispose of an explicit set
$I=\{Q,S(n), n \in \mathbb{Z}\}$ of generating cofibrations and an
explicit set $J=\{ F, R(n), n \in \mathbb{Z}\}$ of generating trivial
cofibrations for $\mathsf{dgcat}$.

Now remark that all conditions of
theorem~\ref{liftarg} are satisfied. In particular
proposition~\ref{path} and the fact that every object in
$\mathsf{dgcat}$ is fibrant, see~\cite{cras}, imply
proposition~\ref{Quillen} which implies that the assumption on cofibrations of
theorem~\ref{liftarg} holds. 

Remark that $F(I)$ is a set of generating
cofibrations on $\mathsf{T}_{\alpha}$-$\mathsf{alg}$ and that $F(J)$ is a set of generating
acyclic cofibrations on $\mathsf{T}_{\alpha}$-$\mathsf{alg}$. This implies that the Quillen
model structure on $\mathsf{T}_{\alpha}$-$\mathsf{alg}$ is cofibrantly
generated. Since the functor $U$ preserves by definition weak
equivalences and fibrations the adjunction $(F,U)$ is a Quillen
adjunction.
This proves the theorem.
\end{proof}

\section{Exact  $\alpha$-cocomplete DG categories}
In this section, we will construct a category
$\mathsf{dgcat}_{ex,\alpha}$ by considering specific diagrams in
$\mathsf{T}_{\alpha}$-$\mathsf{alg}$. The objects of
$\mathsf{dgcat}_{ex,\alpha}$ will be essentially the dg categories
which are stable under suspensions, cosuspensions, cones and
$\alpha$-small sums, see remark~\ref{choices1}.

\begin{definition}\label{exact}
Let $\mathcal{P}$ be the dg category with only one object $X$ and
whose dg algebra of endomorphisms is $k$ concentrated in degree
$0$. Let $\mathcal{S}$, respectively $\mathcal{S}^{-1}$, be the full
sub dg category of $\mathcal{C}_{dg}(\mathcal{P})$, whose objects are
$\widehat{X}$ and $\widehat{X}[1]$, respectively $\widehat{X}$ and
$\widehat{X}[-1]$. We dispose of a fully faithful dg functor $\mathcal{P}
\stackrel{S}{\rightarrow} \mathcal{S}$, respectively $P
\stackrel{S}{\rightarrow} \mathcal{S}^{-1}$.
Let $\mathcal{M}$ be the dg category which has two objects $0$ and $1$
and is generated by a morphism $f$ from $0$ to $1$ that satisfies
$d(f)=0$. Let $\mathcal{C}$ be the full sub dg category of
$\mathcal{C}_{dg}(\mathcal{M})$, whose objects are $\widehat{0}$,
$\widehat{1}$ and $\mathsf{cone}(\widehat{f})$. We dispose of a fully faithful
dg functor $\mathcal{M} \stackrel{C}{\rightarrow} \mathcal{C}$.
\end{definition}
Let $\mathcal{A}$ be a small dg category. 

\begin{remark}\label{choice}
Remark that giving a dg functor $H:\mathcal{S} \rightarrow
\mathcal{A}$, respectively $H':\mathcal{S}^{-1} \rightarrow
\mathcal{A}$, corresponds exactly to specifying two objects $X$ and
$Y$ in $\mathcal{A}$ and an isomorphism $\widehat{X}[1]
  \stackrel{\sim}{\longrightarrow} \widehat{Y}$, respectively $\widehat{X}[-1]
  \stackrel{\sim}{\longrightarrow} \widehat{Y}$, in
  $\mathcal{C}_{dg}(\mathcal{A})$. Remark also that giving a dg
  functor $R:\mathcal{C} \rightarrow \mathcal{A}$ corresponds exactly
  to specifying a morphism $f:X \rightarrow Y$ of degree zero in
  $\mathcal{A}$ such that $d(f)=0$, an object $Z$ in $\mathcal{A}$ and
  an isomorphism $Z \stackrel{\sim}{\longrightarrow} \mathsf{cone}(\widehat{f})$.
\end{remark}

Recall from section~\ref{monad} that we dispose of an adjunction
$$
\xymatrix{
\mathsf{T}_{\alpha}\mbox{-}\mathsf{alg} \ar@<1ex>[d]^U \\
\mathsf{dgcat} \ar@<1ex>[u]^F\,.
}
$$

\begin{definition}\label{exact1}
Let $\mathsf{dgcat}_{ex,\alpha}$ be the category whose objects are the
$4$-tuples $\underline{A}=(A,S_A,S^{-1}_A,C_A)$, where $A$ is a
$\mathsf{T}_{\alpha}$-algebra and $S_A$, $S_A^{-1}$ and $C_A$, the
structure morphisms of $A$, are $\mathsf{T}_{\alpha}$-algebra morphisms
which make both diagrams
$$
\xymatrix{
\underset{\mathcal{P} \rightarrow \mathcal{A}}{\coprod} F(\mathcal{P}) \ar[r]
\ar[d] & A & \underset{\mathcal{P} \rightarrow \mathcal{A}}{\coprod} F(\mathcal{P}) \ar[r]
\ar[d] & A & \underset{\mathcal{M} \rightarrow
  \mathcal{A}}{\coprod} F(\mathcal{M}) \ar[r] \ar[d] & A \\
\underset{\mathcal{P} \rightarrow \mathcal{A}}{\coprod} F(\mathcal{S}) 
\ar[ur]_-{S_A} & & \underset{\mathcal{P} \rightarrow
  \mathcal{A}}{\coprod} F(\mathcal{S}^{-1}) \ar[ur]_-{S^{-1}_A}
& &  \underset{\mathcal{M} \rightarrow \mathcal{A}}{\coprod}
F(\mathcal{C}) \ar[ur]_-{C_A} & 
}
$$
commutative.

A morphism $G: \underline{A} \rightarrow \underline{B}$ in
$\mathsf{dgcat}_{ex,\alpha}$ consists of a morphism of
$\mathsf{T}_{\alpha}$-algebras $G: A \rightarrow B$ that makes the
following diagrams
$$
\xymatrix{
\underset{\mathcal{P} \rightarrow \mathcal{A}}{\coprod} F(\mathcal{S})
\ar[d]_-{S_A} \ar[r] & \underset{\mathcal{P} \rightarrow \mathcal{B}}{\coprod} F(\mathcal{S})
\ar[d]^-{S_B} &  \underset{\mathcal{P} \rightarrow
  \mathcal{A}}{\coprod} F(\mathcal{S}^{-1}) \ar[d]_-{S^{-1}_A} \ar[r]
&   \underset{\mathcal{P} \rightarrow
  \mathcal{B}}{\coprod} F(\mathcal{S}^{-1}) \ar[d]^-{S^{-1}_B}
&  \underset{\mathcal{M} \rightarrow \mathcal{A}}{\coprod}
F(\mathcal{C}) \ar[d]_-{C_A} \ar[r] &  \underset{\mathcal{M} \rightarrow \mathcal{B}}{\coprod}
F(\mathcal{C}) \ar[d]^-{C_B}\\
A \ar[r]_-{G} & B & A \ar[r]_-{G} & B & A \ar[r]_-{G} & B 
}
$$
commutative.
\end{definition}

\begin{remark}\label{choices1}
Remark that an object $\underline{A}$ in $\mathsf{dgcat}_{ex,\alpha}$
consists of a $\mathsf{T}_{\alpha}$-algebra $A$ and of a choice, in 
the sense of remark~\ref{choice}, for the suspensions and cosuspensions
of every object of the dg category $\mathcal{A}=U(A)$ and also of a
choice for the cone of every cycle of degree zero of the dg category
$\mathcal{A}$. In particular $\mathcal{A}$ is stable under
suspensions, cosuspensions and cones. 

Remark also that a morphism $G$
in $\mathsf{dgcat}_{ex,\alpha}$ consists of a morphism of
$\mathsf{T}_{\alpha}$-algebras that commutes with all these choices. 
\end{remark}

We dispose of a forgetful functor
$$ U_1: \mathsf{dgcat}_{ex,\alpha} \longrightarrow
\mathsf{T}_{\alpha}\mbox{-}\mathsf{alg}\,,$$
that associates to an object $\underline{A}$ of
$\mathsf{dgcat}_{ex,\alpha}$ the $\mathsf{T}_{\alpha}$-algebra $A$.

\begin{proposition}
The functor $U_1$ admits a left adjoint.
\end{proposition}

\begin{proof}
The proof will consist in verifying the conditions of theorem~$2$
from section V.$6$. of \cite{Macl}.

We will now prove that the category $\mathsf{dgcat}_{ex,\alpha}$
admits small limits and that these are preserved by the functor
$U_1$. Let $\{\underline{A_i} \}_{i \in I}$ be a family of objects in
$\mathsf{dgcat}_{ex,\alpha}$. Endow the $\mathsf{T}_{\alpha}$-algebra
$\underset{i \in I}{\prod} A_i$ with the structure morphisms induced
by those of $\underline{A_i},\, i \in I$. In this way, the
$\mathsf{T}_{\alpha}$-algebra $\underset{i \in I}{\prod} A_i$ belongs
naturally to $\mathsf{dgcat}_{ex,\alpha}$ and we remark that it is the
product in $\mathsf{dgcat}_{ex,\alpha}$ of the familly
$\{\underline{A_i} \}_{i \in I}$.

Consider now morphisms $G_1,G_2:\underline{A} \rightarrow
\underline{B}$ in $\mathsf{dgcat}_{ex,\alpha}$. Let $K$ be the
equalizer in $\mathsf{T}_{\alpha}$-$\mathsf{alg}$ of the pair
$G_1,G_2:A \rightarrow B$. By remark~\ref{choices1}, we need to show
that the dg category $U(K)$ is endowed with a choice for the
suspension and cosuspension for each object and with a choice for the
cone of every cycle of degree zero. Since $U$ is a right adjoint
functor, $U(K)$ identifies with the equalizer of the pair 
$$U(G_1),U(G_2): \mathcal{A} \rightarrow \mathcal{B}$$
and since $G_1$ and $G_2$
are morphisms in $\mathsf{dgcat}_{ex,\alpha}$, the dg functors
$U(G_1)$ and $U(G_2)$ commute with all the choices. This implies that
the the non-full dg subcategory $U(K)$ of $\mathcal{A}$ is stable
under all the choices of suspension, cosuspension and cones in
$\mathcal{A}$. This shows us that $K$ belongs naturally to
$\mathsf{dgcat}_{ex,\alpha}$ and that it is the equalizer in
$\mathsf{dgcat}_{ex,\alpha}$ of the pair $G_1,G_2:\underline{A}
\rightarrow \underline{B}$.
This proves that the category $\mathsf{dgcat}_{ex,\alpha}$ admits
small limits and that they are preserved by the functor $U$.

We will now prove that the solution set condition is verified,
see~\cite{Macl}.

Let $A$ be a $\mathsf{T}_{\alpha}$-algebra.
Consider the following set of morphisms in
$\mathsf{T}_{\alpha}$-$\mathsf{alg}$
$$ \mbox{Ens}:=\left\{ F(\mathcal{P}) \stackrel{F(S)}{\longrightarrow}
F(\mathcal{S}),\, F(\mathcal{P}) \stackrel{F(S^{-1})}{\longrightarrow}
F(\mathcal{S}^{-1}),\,  F(\mathcal{M}) \stackrel{F(C)}{\rightarrow}
F(\mathcal{C})\right\}\,.$$
Since $\mathcal{P}$ and $\mathcal{M}$ are clearly small in
$\mathsf{dgcat}$, see \cite{cras}, and since the functor $U$
commutes with $\alpha$-filtered colimits, the objects $F(\mathcal{P})$ and
$F(\mathcal{M})$ are $\alpha$-small in
$\mathsf{T}_{\alpha}$-$\mathsf{alg}$, see \cite{Hirschhorn}.

Apply the small object argument, see \cite{Hirschhorn}, to the
morphism
$$ A \longrightarrow 0\,,$$
where $0$ denotes the terminal object in
$\mathsf{T}_{\alpha}$-$\mathsf{alg}$, using the set of morphismes
$\mbox{Ens}$. We dispose of a factorization
$$ 
\xymatrix{
A \ar[rr] \ar[dr]_-i & & 0 \\
 & \mathsf{Ex}_{\alpha}(A) \ar[ur]_-q & \,, 
}
$$
where $q$ is a morphism of $\mathsf{T}_{\alpha}$-algebras that has the
right lifting property with respect to all elements of $\mbox{Ens}$.

Now, for each one of the following (solid) commutative squares
$$
\xymatrix{
F(\mathcal{P}) \ar[r] \ar[d]_{F(S)} & \mathsf{Ex}_{\alpha}(A) \ar[d] &
F(\mathcal{P}) \ar[d]_-{F(S^{-1})} \ar[r] & \mathsf{Ex}_{\alpha}(A)
\ar[d] & F(\mathcal{M}) \ar[d]_-{F(C)} \ar[r] & \mathsf{Ex}_{\alpha}(A)
\ar[d] \\
F(\mathcal{S}) \ar[r] \ar@{.>}[ur] & 0 & F(\mathcal{S}^{-1}) \ar[r]
\ar@{.>}[ur] & 0 & F(\mathcal{C}) \ar@{.>}[ur] \ar[r] & 0 \,,
}
$$
choose a morphism of $\mathsf{T}_{\alpha}$-algebras, (here denoted by a
dashed arrow), as in the proof of proposition $10.5.16$ from
\cite{Hirschhorn}, that makes both triangles commutative. Remark that
a set of morphisms as this one specifies structure
morphisms for the $\mathsf{T}_{\alpha}$-algebra
$\mathsf{Ex}_{\alpha}(A)$. This shows us that when endowed with these
choices $\mathsf{Ex}_{\alpha}(A)$ belongs to $\mathsf{dgcat}_{ex,\alpha}$.

Let now $\underline{B}$ be an object of $\mathsf{dgcat}_{ex,\alpha}$
and $Q:A \rightarrow B$ be a morphism of
$\mathsf{T}_{\alpha}$-algebras.
Remark that the structure morphisms of $B$ and the
construction of the $\mathsf{T}_{\alpha}$-algebra
$\mathsf{Ex}_{\alpha}(A)$ by the small object argument, see the proof of
proposition $10.5.16$ in \cite{Hirschhorn}, allows us to define, by
transfinite induction, a morphism $\overline{Q}$ of
$\mathsf{T}_{\alpha}$-algebras such that the diagram
$$
\xymatrix{
A \ar[r]^-{i} \ar[dr]_-{Q} &  \mathsf{Ex}_{\alpha}(A)
\ar[d]^-{\overline{Q}} \\
 & B
}
$$
commutes. Remark also that when the $\mathsf{T}_{\alpha}$-algebra
$\mathsf{Ex}_{\alpha}(A)$ is endowed with the above choices the
morphism $\overline{Q}$ becomes a morphism in
$\mathsf{dgcat}_{ex,\alpha}$. This proves the solution set condition.

The proposition is now proven.
\end{proof}

We dispose of the following adjunctions
$$
\xymatrix{
\mathsf{dgcat}_{ex,\alpha} \ar@<1ex>[d]^{U_1}\\
\mathsf{T}_{\alpha}\mbox{-}\mathsf{alg} \ar@<1ex>[d]^U
\ar@<1ex>[u]^{F_1} \\
\mathsf{dgcat} \ar@<1ex>[u]^F\,.
}
$$

\begin{proposition}\label{monadic}
The functor $U_1$ is monadic, see \cite{Macl}.
\end{proposition}

\begin{proof}
The proof will consist in verifying condition $(iii)$ of theorem $1$
from section VI.$7$. of \cite{Macl}.

Let $G_1,G_2:\underline{A} \rightarrow \underline{B}$ be a pair of
morphisms in $\mathsf{dgcat}_{ex,\alpha}$. Consider the following
split coequalizer, see \cite{Macl}, in
$\mathsf{T}_{\alpha}$-$\mathsf{alg}$
$$
\xymatrix{
A \ar@<0.5ex>[r]^{G_1} \ar@<-0.5ex>[r]_{G_2} & B \ar@/^1.5pc/[l]^-{R}  \ar[r]^L & D
\ar@/^1.5pc/[l]^-{S} \,,
}
$$
where $L\circ G_1 = L \circ G_2$, $L\circ S =Id$, $G_1 \circ R = Id$,
and $G_2\circ R=S\circ L$.

We will now construct structure morphisms for $D$ such that $D$ will
become an object of $\mathsf{dgcat}_{ex,\alpha}$ and $L$ a morphism
in $\mathsf{dgcat}_{ex,\alpha}$. Apply the functor $U$ to the previous
split coequalizer in $\mathsf{T}_{\alpha}$-$\mathsf{alg}$ and obtain
$$
\xymatrix{
\mathcal{A} \ar@<0.5ex>[r]^{G_1} \ar@<-0.5ex>[r]_{G_2} & \mathcal{B} \ar@/^1.5pc/[l]^-{R}  \ar[r]^L & \mathcal{D}
\ar@/^1.5pc/[l]^-{S} \,.
}
$$
Now, apply the functors
$$ \underset{\mathcal{M} \rightarrow ?}{\coprod} \mathcal{M}\,,
\underset{\mathcal{M} \rightarrow ?}{\coprod}\mathcal{C} :
\mathsf{dgcat} \rightarrow \mathsf{dgcat}$$
to the previous split equalizer in $\mathsf{dgcat}$ and obtain the
following diagram
$$
\xymatrix{
 & \underset{\mathcal{M} \rightarrow \mathcal{A}}{\coprod} \mathcal{C}
 \ar@<0.5ex>[rr] \ar@<-0.5ex>[rr] \ar[ldd]^(.3){C_{\mathcal{A}}}|\hole & & \underset{\mathcal{M} \rightarrow \mathcal{B}}{\coprod} \mathcal{C}
 \ar[rr] \ar[ldd]^(.3){C_{\mathcal{B}}}|\hole & & \underset{\mathcal{M} \rightarrow
   \mathcal{C}}{\coprod} \mathcal{C} \ar@{.>}[ldd]^-{C_{\mathcal{D}}} \\
\underset{\mathcal{M} \rightarrow \mathcal{A}}{\coprod} \mathcal{M}
\ar@<0.5ex>[rr] \ar@<-0.5ex>[rr] \ar[d]  \ar[ur] & & \underset{\mathcal{M} \rightarrow
  \mathcal{B}}{\coprod} \mathcal{M} \ar[rr] \ar[d] \ar[ur] & &
\underset{\mathcal{M} \rightarrow \mathcal{C}}{\coprod} \mathcal{M}
\ar[d]  \ar[ur] & \\
\mathcal{A} \ar@<0.5ex>[rr]^-{G_1}  \ar@<-0.5ex>[rr]_-{G_2} & &
\mathcal{B} \ar[rr]^L  \ar@/^1.5pc/[ll]^-{R}& & \mathcal{D} \ar@/^1.5pc/[ll]^-{S}& 
}
$$

Remark that since $U$ is a split coequalizer, the rows in the
diagram are coequalizers. This implies that the dg functors
$C_{\mathcal{A}}$ and $C_\mathcal{B}$, which correspond under the
adjunction $(F,U)$ to the structure morphisms $C_A$ and $C_B$, induce
a dg functor $C_{\mathcal{D}}$. Now, since $L$ admits a right inverse $S$, a
diagram chasing argument shows us that the right triangle in the
diagram is commutative.
 
Remark that under the adjunction $(F,U)$ this commutative diagram
corresponds exactly to a structure morphism $C_D$ on $D$. Clearly by construction the morphism $U$
preserves this structure morphism. Now, consider an analogous
argument for the construction of structure morphisms $S_D$ and
$S^{-1}_D$.

We will now prove that $L$ is a coequalizer in the category
$\mathsf{dgcat}_{ex,\alpha}$. Let $\underline{E}$ be an object of
$\mathsf{dgcat}_{ex,\alpha}$ and $H:\underline{B} \rightarrow
\underline{E}$ a morphim such that $H \circ G_1 = H \circ G_2$. Since
the morphism $U$ is a coequalizer in
$\mathsf{T}_{\alpha}$-$\mathsf{alg}$, there exists an unique morphism
of $\mathsf{T}_{\alpha}$-algebras $R$ which makes the following
diagram
$$
\xymatrix{
A \ar@<0.5ex>[r]^{G_1} \ar@<-0.5ex>[r]_{G_2} & B \ar[r]^U \ar[dr]_H & D
\ar@{.>}[d]^R \\
 & & E
}
$$
commutative.
We now prove that $R$ is a morphism in
$\mathsf{dgcat}_{ex,\alpha}$. Since $U \circ S = Id$ we have $R=H\circ
S$. Now by apply the functors 
$$\underset{\mathcal{P}
  \rightarrow?}{\coprod}\mathcal{S}\,,\underset{\mathcal{P}
  \rightarrow ?}{\coprod}\mathcal{S}^{-1}\,,\underset{\mathcal{M}
  \rightarrow ?}{\coprod}\mathcal{C} : \mathsf{dgcat} \longrightarrow \mathsf{dgcat}$$
to the image by the functor $U$ of the diagram above and use a diagram
chasing argument to conclude that $U$ belongs in fact to $\mathsf{dgcat}_{ex,\alpha}$. This proves the
proposition.
\end{proof}

Remark that since $U_1$ is monadic and the category $\mathsf{T}_{\alpha}$-$\mathsf{alg}$ is
complete, proposition $4.3.1$ from \cite{Bor} implies that $\mathsf{dgcat}_{ex,\alpha}$ is also complete.

\begin{proposition}\label{filtered}
The category $\mathsf{dgcat}_{ex,\alpha}$ admits $\alpha$-filtered
colimits and these are preserved by the functor $U_1$.
\end{proposition}

\begin{proof}
Let $\{\underline{A_i}\}_{i \in I}$ be a $\alpha$-filtered diagram in
$\mathsf{dgcat}_{ex,\alpha}$. Consider the colimit
$$ Y :=\underset{i \in I}{\mbox{colim}}\,A_i\,,$$
of the $\alpha$-filtered diagram $\{A_i\}_{i \in I}$ of
$\mathsf{T}_{\alpha}$-algebras. 

We will now construct structure
morphisms for $Y$ such that $Y$ becomes the colimit in
$\mathsf{dgact}_{ex,\alpha}$ of the diagram $\{\underline{A_i}\}_{i
  \in I}$. Since the functor $U$ commutes with $\alpha$-filtered
colimits, see section~\ref{homotopy}, we have
$$ U(Y) = \underset{i \in I}{\mbox{colim}}\, \mathcal{A}_i\,.$$
We now construct a structure morphism $C_Y$, see
definition~\ref{exact1}, for $Y$. We dispose of the following
$\alpha$-filtered diagrams in $\mathsf{dgcat}$ and morphismes
between them
$$
\xymatrix{
\left\{ \underset{\mathcal{M} \rightarrow
  \mathcal{A}_i}{\coprod}\mathcal{M} \right\}_{i \in I} \ar[r] \ar[d] & \left\{
\mathcal{A}_i \right\}_{i \in I} \\
\left\{ \underset{\mathcal{M} \rightarrow
  \mathcal{A}_i}{\coprod}\mathcal{C} \right\}_{i \in I}
\ar[ur]_-{C_{\mathcal{A}_i}} & 
}
$$
Now remark that since we are considering $\alpha$-filtered colimits,
we have
$$ \underset{i \in I}{\mbox{colim}} \underset{\mathcal{M}\rightarrow
  \mathcal{A}_i}{\coprod}\mathcal{M} \stackrel{\sim}{\longrightarrow} \underset{\mathcal{M}
  \rightarrow \underset{i \in I}{\mbox{colim}}\,\mathcal{A}_i}{\coprod}
\mathcal{M}$$
and
$$ \underset{i \in I}{\mbox{colim}} \underset{\mathcal{M}\rightarrow
  \mathcal{A}_i}{\coprod}\mathcal{C} \stackrel{\sim}{\longrightarrow} \underset{\mathcal{M}
  \rightarrow \underset{i \in I}{\mbox{colim}}\,\mathcal{A}_i}{\coprod}
\mathcal{C}\,.$$
This implies that the morphism of $\mathsf{T}_{\alpha}$-algebras which
corresponds under the adjunction $(F,U)$ to the dg functor
$ \underset{i \in I}{\mbox{colim}} \,C_{\mathcal{A}_i} $
is a structure morphism $C_Y$ of $Y$. Consider now an analogous
argument for the construction of structure morphisms $S_Y$ and
$S^{-1}_Y$.

Finally since the functor $U$ commutes with
$\alpha$-filtered colimits, $Y$ is clearly the colimit in
$\mathsf{dgcat}_{ex,\alpha}$ of the diagram  $\{\underline{A_i}\}_{i
  \in I}$. This proves the proposition.
\end{proof}

\begin{proposition}
The category $\mathsf{dgcat}_{ex,\alpha}$ is cocomplete.
\end{proposition}

\begin{proof}
Recall that by proposition~\ref{monadic} the adjunction
$$
\xymatrix{
\mathsf{dgcat}_{ex,\alpha} \ar@<1ex>[d]^{U_1}\\
\mathsf{T}_{\alpha}\mbox{-}\mathsf{alg}
\ar@<1ex>[u]^{F_1} \,,
}
$$
is monadic. Now by propositions $4.3.2$ and $4.3.6$ of \cite{Bor} we
only need to show that the functor $U_1 \circ F_1$ commutes with
$\alpha$-filtered colimits. But this follows from the fact that $F$ is
a left adjoint and that by proposition~\ref{filtered} $\alpha$-filtered
colimits exist in $\mathsf{dgcat}_{ex,\alpha}$ and are preserved by
$U_1$.
This proves the proposition.
\end{proof}

From now on and until the end of this section we consider definition~\ref{Quillen} applied to our particular adjunction
$(F_1,U_1)$, see theorem~\ref{main}.

Let $\underline{A}$ be an object of $\mathsf{dgcat}_{ex,\alpha}$.
\begin{proposition}\label{path1}
The $\mathsf{T}_{\alpha}$-algebra $P(A)$, see proposition~\ref{path},
is endowed with natural structure morphisms and so $\underline{A}$
admits a path object in $\mathsf{dgcat}_{ex,\alpha}$.
\end{proposition}
\begin{proof}
We will construct structure morphisms $S_{P(A)}$, $S_{P(A)}^{-1}$ and
$C_{P(A)}$ for $P(A)$, see definition~\ref{exact1}, in such a way that
$P(A)$ becomes a path object in $\mathsf{dgcat}_{ex,\alpha}$.
Construct the structure morphism $S_{P(A)}$, respectively
$S_{P(A)}^{-1}$, by applying $S_A$, respectively $S_A^{-1}$, in each
component of the objects of $P(A)$.

By remark~\ref{choices1}, to construct a
structure morphism $C_{P(A)}$, we need to construct a cone in the
category $P(\mathcal{A})$ for every cycle of degree zero in $P(\mathcal{A})$. Let now $(m_X,m_Y,h)$ be a cycle of degree $0$ in $P(\mathcal{A})$
between the objects $X \stackrel{f}{\rightarrow} Y$ and $X'
\stackrel{f'}{\rightarrow} Y'$ of $P(\mathcal{A})$. In particular $h$
is a morphism in $\mathcal{A}$ of degree $-1$ such that $d(h)=m_Y\circ
f - f'\circ m_X$. Consider the following diagram
$$
\xymatrix{
X \ar[d]_{m_X} \ar[rr]^f  \ar[drr]|{h}& &  Y \ar[d]^{m_Y} \\
X' \ar[rr]^{f'} \ar[d]_{i_X}  \ar@{.>}[drr]|{h'} &  & Y' \ar[d]^{i_Y} \\
\mathsf{cone}(m_X) \ar@{.>}[rr]_-{\Phi} &  & \mathsf{cone}(m_Y)\,,
}
$$
where $h'=0$ and $\Phi$ is the morphism defined by the matrix
$$
\begin{bmatrix} 
f' & h \\ 0 & S(f)
\end{bmatrix}\,.
$$
Remark that the object in $P(\mathcal{A})$
$$ 
\xymatrix{
\mathsf{cone}(m_X) \ar@{.>}[rr]^-{\Phi} & & \mathsf{cone}(m_Y)\,,
}
$$
corepresents the cone of the morphism $(m_X,m_Y,h)$ in $P(\mathcal{A})$.
This shows us that $P(A)$ belongs to $\mathsf{dgcat}_{ex,\alpha}$. It
is also clear that the morphisms $A
\stackrel{i_A}{\longrightarrow} P(A)$ and 
$$
\xymatrix{
 P(A) \ar@{->>}[r]^-{q_A} &
 A\times A
}
$$ are in fact morphisms in
$\mathsf{dgcat}_{ex,\alpha}$. This proves the proposition.
\end{proof}

\begin{theorem}
The category $\mathsf{dgcat}_{ex,\alpha}$ when endowed with
the notions of weak equivalence, fibration and cofibration as in
definition~\ref{lift}, becomes a cofibrantly generated Quillen model
category and the adjunction $(F_1,U_1)$ becomes a Quillen adjunction.
\end{theorem}

\begin{proof}
Recall from theorem~\ref{main} that we dispose of an explicit set
$F(I)$ of generating cofibrations and an
explicit set $F(J)$ of generating trivial
cofibrations for $\mathsf{T}_{\alpha}$-$\mathsf{alg}$. 

Now remark that all conditions of
theorem~\ref{liftarg} are satisfied. In particular
proposition~\ref{path1} and the fact that every object in
$\mathsf{T}_{\alpha}$-$\mathsf{alg}$ is fibrant imply
proposition~\ref{Quillen} which implies that the assumption on cofibrations of
theorem~\ref{liftarg} holds. 

Remark that $F_1(F(I))$ is a set of generating
cofibrations on $\mathsf{dgcat}_{ex,\alpha}$ and that $F_1(F(J))$ is a set of generating
acyclic cofibrations on $\mathsf{dgcat}_{ex,\alpha}$. This implies that the Quillen
model structure on $\mathsf{dgcat}_{ex,\alpha}$ is cofibrantly
generated. Since the functor $U_1$ preserves by definition weak
equivalences and fibrations the adjunction $(F_1,U_1)$ is a Quillen
adjunction.
This proves the theorem.
\end{proof}

\section{Enhancement of well-generated algebraic triangulated
  categories}

In this section, we show that the category of $\alpha$-compactly generated algebraic
triangulated categories up to equivalence admits a natural Quillen
enhancement given by our model category $\mathsf{dgcat}_{ex,\alpha}$.

Let $\mathsf{Tri}_{\alpha}$ denote the category of $\alpha$-compactly generated algebraic
triangulated categories in the sense of Neeman, see \cite{Krause} \cite{Neeman}.

Let $\underline{A}$ be an object of $\mathsf{dgcat}_{ex,\alpha}$, with
underlying dg category $\mathcal{A}$.
Recall from proposition~\ref{sums} that $\mathcal{A}$ admits
$\alpha$-small sums.

\begin{definition}[\cite{Porta}]
The $\alpha$-continuous derived category $\mathcal{D}_{\alpha}(\underline{A})$ of
$\underline{A}$ is the triangle quotient of the derived category
$\mathcal{D}(\mathcal{A})$ of $\mathcal{A}$ by the localizing
triangulated subcategory of $\mathcal{D}(\mathcal{A})$ generated by
the cones on the canonical morphisms
$$ \underset{i \in I}{\bigoplus} \widehat{X_i} \longrightarrow
\widehat{\underset{i \in I}{\bigoplus} X_i}\,,$$
where $(X_i)_{i \in I}$ is a family of objects of $\mathcal{A}$ and $I$ is a set of
cardinality strictly smaller than $\alpha$.
\end{definition}

\begin{remark}
By a theorem of \cite{Porta}, $\mathcal{D}_{\alpha}(\underline{A})$ belongs to
$\mathsf{Tri}_{\alpha}$. 

Remark that this construction is functorial in
$A$. In fact, let $F:\underline{A} \rightarrow \underline{B}$ be a
morphism in $\mathsf{dgcat}_{ex,\alpha}$. Since the dg functor
$F:\mathcal{A} \rightarrow \mathcal{B}$ commutes with $\alpha$-small
sums, see remark~\ref{sums2}, it induces a functor
$$ \mathcal{D}_{\alpha}(F): \mathcal{D}_{\alpha}(\underline{A})
\rightarrow D_{\alpha}(\underline{B})$$
between triangulated categories.
\end{remark}
Thus, we have defined a functor
$$ \mathcal{D}_{\alpha}(-): \mathsf{dgcat}_{ex,\alpha}
\longrightarrow \mathsf{Tri}_{\alpha}\,.$$
Remark that if $F:\underline{A} \rightarrow \underline{B}$ is a weak
equivalence in $\mathsf{dgcat}_{ex,\alpha}$, i.e. the dg functor
$F:\mathcal{A} \rightarrow \mathcal{B}$ induces an equivalence of
categories $\mathsf{H}^0:\mathsf{H}^0(\mathcal{A}) \rightarrow
\mathsf{H}^0(\mathcal{B})$, then the triangulated functor
$\mathcal{D}_{\alpha}(F)$ is an equivalence of triangulated categories.

The following theorem is proven in \cite{Porta}, \emph{cf.}~\cite{Ober}.

\begin{theorem}[\cite{Porta}]
The functor $\mathcal{D}_{\alpha}(-)$ satisfies the conditions:
\begin{itemize}
\item[-]  every category $\mathcal{T}$ in $\mathsf{Tri}_{\alpha}$ is equivalent to
  $\mathcal{D}_{\alpha}(\underline{A})$ for some object
  $\underline{A}$ in $\mathsf{dgcat}_{ex,\alpha}$ and
\item[-] a morphism $F$ in $\mathsf{dgcat}_{ex,\alpha}$ is a weak
  equivalence if and only if $\mathcal{D}_{\alpha}(F)$ is an
  equivalence of triangulated categories.
\end{itemize}
\end{theorem}

\end{document}